\documentclass[twoside]{article}
\usepackage{style}

\author{Mireille Boutin and Gregor Kemper}
\title{Which Point Configurations are Determined by the 
Distribution of their Pairwise Distances?}
\begin{document}
\maketitle

\begin{abstract}
  In a previous paper we showed that, for any $n \ge m+2$, most sets
  of~$n$ points in $\RR^m$ are determined (up to rotations,
  reflections, translations and relabeling of the points) by the
  distribution of their pairwise distances. But there are some
  exceptional point configurations which are not reconstructible from
  the distribution of distances in the above sense. In this paper,
  we present a reconstructibility test with running time $O(n^{11})$.
  The cases of
  orientation preserving rigid motions (rotations and translations)
  and scalings are also discussed.
\end{abstract}

\section*{Introduction} \label{0s:myintro}
In this paper, we present a quick and easy 
(but slightly imperfect) 
solution to the problem of
characterizing the shape of sets of $n$ points in Euclidean space, 
so-called $n$-point configurations,
for any positive integer $n$.
More precisely, an {\em $n$-point configuration}
is a collection of $n$ points in $\RR^m$.
Point configurations often arise in biological and medical imagery,
as well as in the fields of archaeology, astronomy and cartography, 
to name just a few.
For example,
stellar constellations, minutiae of fingerprints, 
and distinguished points (landmarks) on medical images
represent point configurations.

An important problem of computer vision is that of 
recognizing point configurations.
In other words, the problem is to determine whether
two point configurations have the same shape,
that is to say, whether there exists a rotation and a translation 
(sometimes a reflection and/or a scaling are allowed as well)
which maps the first point configuration onto the second.
Let us first concentrate on the case of rigid motions,
i.e.~rotations, translations and reflections in $\RR^m$.
Note that any rigid motion can be written as $(M,T)$,
where $M$ is an orthogonal $m$-by-$m$ matrix
and $T$ is an $m$-dimensional (column) vector.

One of the biggest difficulties in trying to identify point configurations
up to rigid motions
is the absence of labels for the points: 
one does not know, a priori, which point is going to be mapped to which.
If the points were already labeled in correspondence, 
then, following the so-called Procrustes approach (\mycite{Gower:1975}), 
one could analytically determine 
a rigid motion
which maps the first string as close as possible 
(in the $L^2$ sense, for example) to the second. 
The statistical analysis of such methods is presented in \mycite{Goodall}.
Another way to proceed would be to 
compare the pairwise (labeled) distances between the points
of each point configurations (\mycite{Blumenthal}).
Indeed, the following well known fact holds.
See, for example, \mycite{FirstPaper} for a simple proof.
\begin{prop}
\label{distance_geometry}
Let $p_1,\ldots,p_n$ and $q_1,\ldots,q_n$ be points in $\RR^m$.
If $\|p_i-p_j \|=\|q_i-q_j \|$ 
for every $i,j=1,\ldots,n$, 
then there exists a rigid motion $(M,T)$
such that $M p_i+T=q_i$, for every $i=1,\ldots,n$.
\end{prop}

A variety of methods have been developed for labeling the points 
of two $n$-point configuration in correspondence.
See, for example, \mycite{HartleyZisserman} 
for a description of some of these methods.
But labeling the points is a complex task
which we would much rather do without.
Invariant theory 
suggests a possible approach for recognizing unlabeled points.
The idea consists in comparing certain functions 
of the pairwise distances between the points of the configuration
which have the property that they 
are unchanged by a relabeling of the points.
These are often called {\em graph invariants} and 
have been computed in the case of $n=4$ by \mycite{ACG}, 
and $n=5$ by the second author~[\citenumber{Derksen:Kemper}, page~220].
Unfortunately, 
the case $n=6$ or larger still stands as a computational challenge.
Moreover, 
the invariants used are polynomial functions of the distances
whose number and degrees increase dramatically with $n$.
They are thus very sensitive to round off errors and noise.

In the following, 
we study an alternative approach 
based on the use of a very simple object: 
the distribution of the pairwise distances. 
The distribution of the pairwise distances  
of an $n$-point configuration is 
an array which lists all the different values 
of the pairwise distances between the points
in increasing order 
and the number of times each value occurs.
For example, 
the distribution of distances of
four points situated at the corners of a unit square is given in Table 1.

\begin{table}
\centering          
\begin{tabular}{|cc|}      
\hline  
  Value & \# of Occurrences             \\    
\hline  
  $1$    &           $4$   \\ 
 $\sqrt{2}$  &            $2$   \\
 \hline
\end{tabular}
\label{TableSquare}          
\caption{Distribution of distances of a unit square.} 
\end{table}

Obviously, such a distribution remains unchanged under
any rigid motion of the point configuration
as well as any relabeling of the points.
For $n=1,2$ or $3$, it is easy to see that
the distribution of distances completely characterize
the $n$-point configuration up to a rigid motion. 
For $n\geq m+2$,
we proved that, {\em most of the time},
this distribution {\em completely} characterizes the shape of the 
point configuration (see~[\citenumber{FirstPaper}, Theorem~2.6]).  

To simplify our discussion, 
we introduce the concept of reconstructibility from distances.
\begin{defi}
We say that the $n$-point configuration
represented by $p_1,\ldots,p_n\in\RR^m$
is {\em reconstructible from distances} 
if, for every $q_1,\ldots,q_n\in \RR^m$ having 
the same distribution of distances,
there exists a rigid motion $(M,T)$ 
and a permutation $\pi$ of the labels $ \{1,\ldots,n \} $
such that $M p_i+T=q_{\pi (i)}$, for every $i=1,\ldots,n$.
\end{defi}
In the following, 
we shall often identify a point configuration 
and one of its representation $p_1,\ldots,p_n\in\RR^m$.
This is done for simplicity and we hope it will not create any confusion.
Please note that the question of reconstructing the point configuration
from its distribution of distances will 
{\em not} be addressed in this paper. 
We suspect this is quite a challenging problem. In fact, our guess is
that this problem lies in the complexity class NP; it might even be
NP-complete.

Theorem~2.6 of \mycite{FirstPaper} actually implies that 
there exists an open and dense subset $\Omega\subset(\RR^m)^n$
of reconstructible point configurations.
In Section \ref{sReconstructible}, we concentrate on the planar case $m=2$
and give an algorithm  in $O(n^{11})$ steps 
to determine whether a point lies in $\Omega$. 
(A simple Matlab implementation of this algorithm is given in the appendix.)
A generalization to other dimensions $m$ is also mentioned.
Section 3 describes how an additional distribution 
can be used in the planar case
in order to compare
the orientation of two point configurations.
In Section 4, 
we show that a slightly modified distribution 
can be used to completely characterize most point configurations
up to rigid motions and scalings.

\section{Reconstructible configurations}
\label{sReconstructible}

Denote by $\mathcal{P}$ the set of pairs
\[
\mathcal{P}=\left\{ \{i,j \} | i\neq j, i,j=1,\ldots,n \right\}
\]
Consider the group of permutations $S_{\binom{n}{2}}$
of the elements of  $\mathcal{P}$.
For any $\varphi\in S_{\binom{n}{2}}$ 
and any $\{i,j \}\in {\mathcal P}$, we denote by 
$\varphi \cdot \{i,j \}$ 
the image of $\{i,j\}$ under $\varphi$.
For two point configurations to have the same distribution of distances
means that there exists a permutation $\varphi\in S_{\binom{n}{2}}$
which maps the labeled pairwise distances of the first configuration 
onto the labeled pairwise distances of the second configuration.
More precisely, if $p_1,\ldots,p_n\in\RR^m$ and $q_1,\ldots,q_n\in\RR^m$ 
have the same distribution of distances,
let $d_{i,j}=\| p_i-p_j \|^2$ 
and $d_{i,j}'=\| q_i-q_j \|^2$, for all $\{ i,j\}\in {\mathcal P}$. 
Then there exists $\varphi\in S_{\binom{n}{2}}$ such that
\[
d_{\varphi\cdot\{i, j\}}=d_{\{ i,j\}}', 
\text{ for all }\{i,j \}\in {\mathcal P}.
\]

For close enough point configurations, 
we have proved in [\citenumber{FirstPaper}] that
one does not need to keep track of the labeling of the points.
The proof is very short and we reproduce it here for completeness.
\begin{prop}
For any $n$-point configuration $p_1,\ldots,p_n\in\RR^m$,
there exists a neighborhood $U$ of \newline $(p_1,\ldots,p_n)\in(\RR^m)^n$
such that if
$(q_1,\ldots,q_n)\in U$ is an $n$-point configuration
with the same distribution of distances as that of 
$(p_1,\ldots,p_n)$, then the two point configurations
are the same up to a rigid motion and a relabeling
of the points.
\end{prop}
\begin{proof}
 Let us assume the contrary. Then
 there exists a sequence of $n$-point configurations 
 $\{ q_1^k,\ldots,q_n^k \}_{k=1}^\infty $
 converging to $p_1,\ldots,p_n$,
 and a sequence of permutations 
$\{ \varphi_k \}_{k=1}^\infty\subset S_{\binom{n}{2}}$ 
 such that none of the
 $q_1^k,\ldots,q_n^k$ 
 can be mapped to $p_1,\ldots,p_n$ 
 by a rotation and a translation and a relabeling,
 but the distances 
 $d_{\{i,j \}}=\| p_i-p_j \|^2$
 are mapped to
 the distances
 $d_{\{i,j \}}^{k}=\| q_i^k-q_j^k \|^2$
 by $\varphi_k$ 
 so   $d_{\varphi_k \cdot \{i,j \}}=d_{ \{i,j \}}^{k}$ 
 for all $\{ i,j\} \in {\mathcal P}$.
 By taking a subsequence, we may assume that
 $\varphi_k=\varphi$ is the same for every $k$
 since $S_{\binom{n}{2}}$ is a finite group.
 Taking the limit, we obtain that $ d_{\varphi \cdot \{i,j \}} = 
 \lim_{k\rightarrow \infty} d_{\{i,j \}}^{k}$, 
 for  $\{ i,j\}\in {\mathcal P}$.
 By continuity of the distance,
 this implies that $ d_{\varphi \cdot \{i,j \}} = d_{\{i,j \}}$,
 for all $\{i,j \}\in {\mathcal P}$.
 Therefore, 
 $d_{ \{i,j \}}=d_{ \{i,j \}}^{k}$
 for every $\{ i,j \} \in {\mathcal P}$.
 By Proposition \ref{distance_geometry},
 this implies that
 $q_1^k,\ldots,q_n^k$ and $p_1,\ldots,p_n$
 are the same up to a rigid motion, for every $k$,
 which contradicts our hypothesis, and the conclusion follows.
\end{proof}

Unfortunately, the size of the neighborhood is unknown and varies with
the points $p_1,\ldots,p_n$, so this local result is not very
practical.  We now consider the global case.  Observe that some of the
permutations in $S_{\binom{n}{2}}$ correspond to a relabeling of the
points.  More precisely, $\varphi$ corresponds to a relabeling of the
points if there exists a permutation $\pi:\{1,\ldots,n \}
\hookrightarrow \{1,\ldots,n\}$ of the indices such that $\varphi
\cdot \{ i,j\} = \{ \pi(i),\pi(j)\} $, for every $\{ i,j \} \in
{\mathcal P} $.  Relabelings are the {\em good } permutations: if the
permutation mapping the labeled pairwise distances of a point
configuration onto the labeled pairwise distances of another
configuration is a relabeling, then the two configurations are the
same up to a rigid motion.  We need to know what distinguishes the
good permutations from the bad permutations. The following lemma,
which is central to our argument, says that, informally speaking, a
permutation is a relabeling if it preserves adjacency.

\begin{lemma} \label{lAdjacency}
  Suppose $n \ne 4$. A permutation $\phi \in S_{\binom{n}{2}}$ is a
  relabeling if and only if for all pairwise distinct indices $i,j,k
  \in \{1 \upto n\}$ we have
  \begin{equation} \label{eqAdjacency}
    \phi \cdot \{i,j\} \cap \phi \cdot \{i,k\} \ne \emptyset.
  \end{equation}
\end{lemma}

\begin{proof}
  For $n \le 3$, every $\phi \in S_{\binom{n}{2}}$ is a relabeling,
  and the condition~\eqref{eqAdjacency} is always satisfied. Thus we
  may assume $n \ge 5$. It is also clear that every relabeling
  satisfies~\eqref{eqAdjacency}.
  
  Suppose that $\phi \in S_{\binom{n}{2}}$ is a permutation of
  $\mathcal P$ which satisfies~\eqref{eqAdjacency}. Take any $i,j,k,l
  \in \{1 \upto n\}$ pairwise distinct and assume, by way of
  contradiction, that $\phi \cdot \{i,j\} \cap \phi \cdot \{i,k\} \cap
  \phi \cdot \{i,l\} = \emptyset$. Then the injectivity of $\phi$ and
  the condition~\eqref{eqAdjacency} imply that we can write $\phi
  \cdot \{i,j\} = \{a,b\}$, $\phi \cdot \{i,k\} = \{a,c\}$, and $\phi
  \cdot \{i,l\} = \{b,c\}$ with $a,b,c \in \{1 \upto n\}$ pairwise
  distinct. Now choose $m \in \{1 \upto n\} \setminus \{i,j,k,l\}$.
  Then $\phi \cdot \{i,m\}$ must meet each of the sets $\{a,b\}$,
  $\{a,c\}$, and $\{b,c\}$. Being itself a set of two elements, $\phi
  \cdot \{i,m\}$ must be one of the sets $\{a,b\}$, $\{a,c\}$, or
  $\{b,c\}$, contradicting the injectivity of $\phi$. Therefore $\phi
  \cdot \{i,j\} \cap \phi \cdot \{i,k\} \cap \phi \cdot \{i,l\} \ne
  \emptyset$.
  
  Fix an index $i \in \{1 \upto n\}$ and choose $j,k \in \{1 \upto n\}
  \setminus \{i\}$. Then $\phi \cdot \{i,j\} \cap \phi \cdot \{i,k\}$
  is a set with one element, and by the above this one element must
  also lie in every $\phi \cdot \{i,l\}$ with $l \in \{1 \upto n\}
  \setminus \{i\}$. Hence $\bigcap_{l \ne i} \phi \cdot \{i,l\} \ne
  \emptyset$. This allows us to define a map $\map{\sigma}{\{1 \upto
    n\}}{\{1 \upto n\}}$ with
  \begin{equation} \label{defsigma}
    \bigcap
    \begin{Sb}
      j=1 \\
      j \ne i
    \end{Sb}
    ^n \phi \cdot \{i,j\} = \{\sigma(i)\}.
  \end{equation}
  For $i \in \{1 \upto n\}$ define $M_i := \left\{\{i,j\} \mid j \in
    \{1 \upto n\} \setminus \{i\}\right\}$. Then~\eqref{defsigma}
  tells us that $\phi \cdot M_i \subseteq M_{\sigma(i)}$. Since $|M_i|
  = |M_{\sigma(i)}|$ and since $\phi$ is injective, this implies $\phi
  \cdot M_i = M_{\sigma(i)}$. Take $i,i' \in \{1 \upto n\}$ with
  $\sigma(i) = \sigma(i')$. Then $\phi \cdot M_i = \phi \cdot M_{i'}$,
  which implies $M_i = M_{i'}$ and therefore $i = i'$. Thus $\sigma$
  is injective.
  
  Equation~\eqref{defsigma} implies that for $i,j \in \{1 \upto n\}$
  distinct we can write $\phi \cdot \{i,j\} =
  \{\sigma(i),\gamma_i(j)\}$ with $\map{\gamma_i}{\{1,\ldots ,n\}
    \setminus \{i\}}{\{1,\ldots,n \}}$. But applying~\eqref{defsigma}
  with the roles of~$i$ and~$j$ interchanged yields
  \[
  \{\sigma(j)\} = \bigcap
  \begin{Sb}
    i=1 \\
    i \ne j
  \end{Sb}
  ^n \phi \cdot \{i,j\} = \bigcap
  \begin{Sb}
    i=1 \\
    i \ne j
  \end{Sb}
  ^n \{\sigma(i),\gamma_i(j)\}.
  \]
  By the injectivity of $\sigma$ this implies $\sigma(j) =
  \gamma_i(j)$ for all $i \ne j$. We conclude that $\phi \cdot \{i,j\}
  = \{\sigma (i),\sigma(j)\}$ for all $i,j \in \{1 \upto n\}$
  distinct. But this means that $\phi$ is a relabeling, as claimed
\end{proof}

\begin{rem*}
  For $n = 4$, \lref{lAdjacency} becomes false. An example is given by
  $\phi \in S_{\binom{4}{2}}$ defined as
  \[
  \begin{array}{lll}
    \phi \cdot \{1,2\} = \{1,2\}, & \phi \cdot \{1,3\} = \{1,3\}, &
    \phi \cdot \{1,4\} = \{2,3\} \\
    \phi \cdot \{2,3\} = \{1,4\}, & \phi \cdot \{2,4\} = \{2,4\}, &
    \phi \cdot \{3,4\} = \{3,4\}.
  \end{array}
  \]
  This permutation satisfies~\eqref{eqAdjacency}, but it is not a
  relabeling. \lref{lAdjacency} becomes true for $n = 4$ if we add the
  additional condition
  \begin{equation}
\label{eqTest2}
  \phi \cdot \{1,2\} \cap \phi \cdot \{1,3\} \cap \phi \cdot \{1,4\}
  \ne \emptyset.
  \end{equation}
\end{rem*}

Do non-reconstructible point configurations exist? The answer is yes.
Some examples can be found in \mycite{FirstPaper}.  Fortunately,
non-reconstructible configurations are rare.  The key to this fact is
contained in the functional relationships between the pairwise
distances of a point configuration. These relationships are well-known
from classical invariant theory. For example, a planar configuration
of four points $p_i$, $p_j$, $p_k$, and $p_l$ satisfies
\[
\det\left(
  \begin{array}{cccc}
    -2 d_{\{i,l\}} & d_{\{i,j\}} - d_{\{i,l\}} - d_{\{j,l\}} &
    d_{\{i,k\}} - d_{\{i,l\}} - d_{\{k,l\}} \\
    d_{\{i,j\}} - d_{\{i,j\}} - d_{\{j,l\}} & -2 d_{\{j,l\}} &
    d_{\{j,k\}} - d_{\{j,l\}} - d_{\{k,l\}} \\
    d_{\{i,k\}} - d_{\{i,l\}} - d_{\{k,l\}} & d_{\{j,k\}} -
    d_{\{j,l\}} - d_{\{k,l\}} & -2 d_{\{k,l\}}
  \end{array}
  \right) = 0.
\]
We can also express this relationship as follows. Define the polynomial
\begin{multline*}
  g(U,V,W,X,Y,Z) := 2 U^2 Z + 2 U V X - 2 U V Y - 2 U V Z - 2 U X W -
  2 U X Z + 2 U Y W - \\
  2 U Y Z -  2 U W Z + 2 U Z^2 + 2 V^2 Y - 2 V X Y - 2 V X W + 2 V Y^2
  - 2 V Y W - \\
  2 V Y Z + 2 V W Z + 2 X^2 W - 2 X Y W + 2 X Y Z + 2 X W^2 - 2 X W Z.
\end{multline*}
Then
\begin{equation} \label{eqRelation}
  g\left(d_{\{i,j\}},d_{\{i,k\}},d_{\{i,l\}},d_{\{j,k\}},d_{\{j,l\}},
    d_{\{k,l\}}\right) = 0.
\end{equation}

For simplicity, we continue to concentrate on the planar case $m=2$
although other dimensions can be treated similarly. 
Recall that $\mathcal{P}$ denotes the set of pairs
$
\mathcal{P}=\left\{ \{i,j \} | i\neq j, i,j=1,\ldots,n \right\}
$.
The following theorem gives a practical test for reconstructibility of
planar point configurations.

\begin{theorem} \label{test}
  Let $n\geq 5$, let $p_1 \upto p_n\in \RR^2$ and let $d_{\{ i,j \}
  }=\| p_i-p_j \|^2$ be the square of the Euclidean distance between
  $p_i$ and $p_j$, for every $\{ i,j \}\in {\mathcal P}$.  Suppose
  that for each choice of indices
  $i_0$,$i_1$,$i_2$,$j_1$,$j_2$,$k_1$,$k_2$,$l_1$,\newline $l_2$,$m_1$,$m_2 \in \{1 \upto n\}$ such
  that the pairs $\{i_0,i_1 \}$, $\{i_0,i_2\}$, $\{j_1,j_2\}$,
  $\{k_1,k_2\}$, $\{l_1,l_2\}$, $\{m_1,m_2\} \in \mathcal{P}$ are
  distinct, we have
  \begin{equation} \label{eqTest}
    g\left(d_{\{i_0,i_1\}},d_{\{j_1,j_2\} },d_{\{ k_1,k_2\}},
    d_{\{l_1,l_2\}},d_{\{m_1,m_2\}},d_{\{i_0,i_2\}}\right)\ne 0.
  \end{equation}
  Then $p_1 \upto p_n$ is reconstructible from distances.
\end{theorem}

\begin{proof}
  Let $q_1 \upto q_n \in \RR^2$ be a point configuration with the same
  distribution of distances as $p_1 \upto p_n$. Write
  $d_{\{i,j\}}^\prime = \| q_i-q_j \|^2$. Then there exists a
  permutation $\phi \in S_{\binom{n}{2}}$ of the set $\mathcal P$ such
  that
  \[
  d_{\{i,j\}}^\prime = d_{\phi \cdot \{i,j\}}.
  \]
  We wish to use \lref{lAdjacency} for showing that $\phi^{-1}$ is a
  relabeling, which will imply that $\phi$ is also a relabeling. Take
  any pairwise distinct indices $i,j,k,l \in \{1 \upto n\}$. Then the
  above equation and~\eqref{eqRelation} imply
  \begin{multline*}
    g\left(d_{\phi \cdot \{i,j\}},d_{\phi \cdot \{i,k\}}, d_{\phi
        \cdot \{i,l\}},d_{\phi \cdot \{j,k\}},d_{\phi \cdot
        \{j,l\}},d_{\phi \cdot \{k,l\}}\right) = \\
    g\left(d_{\{i,j\}}^\prime,
      d_{\{i,k\}}^\prime,d_{\{i,l\}}^\prime,d_{\{j,k\}}^\prime,
      d_{\{j,l\}}^\prime,d_{\{k,l\}}^\prime\right) = 0.
  \end{multline*}
  It follows from the hypothesis~\eqref{eqTest} that $\phi \cdot
  \{i,j\}$ and $\phi \cdot \{k,l\}$ are disjoint (otherwise they would
  have an index $i_0$ in common). So for disjoint sets $\{i,j\}$ and
  $\{k,l\}$ we have that $\phi \cdot \{i,j\}$ and $\phi \cdot \{k,l\}$
  are also disjoint. This is equivalent to saying that if $\phi \cdot
  \{i,j\}$ and $\phi \cdot \{k,l\}$ have non-empty intersection, then
  the same is true for $\{i,j\}$ and $\{k,l\}$. Take $a,b,c \in \{1
  \upto n\}$ pairwise distinct and set $\{i,j\} := \phi^{-1} \cdot
  \{a,b\}$ and $\{j,k\} := \phi^{-1} \cdot \{a,c\}$. Then $\phi \cdot
  \{i,j\} \cap \phi \cdot \{k,l\} = \{a,b\} \cap \{a,c\} = \{a\}$,
  hence, as seen above, $\{i,j\}$ and $\{k,l\}$ have non-empty
  intersection. Thus the condition~\eqref{eqAdjacency} of
  \lref{lAdjacency} is satisfied for $\phi^{-1}$. It follows that
  $\phi^{-1}$, and hence also $\phi$, is a relabeling: $\phi \cdot
  \{i,j\} = \{\pi(i),\pi(j)\}$ with $\pi \in S_n$. Now it follows from
  \pref{distance_geometry} that there exists a rigid motion $(M,T)$
  such that
  \[
  q_{\pi(i)} = M p_i + T
  \]
  for all $i \in \{1 \upto n\}$. This completes the proof.
\end{proof}

\begin{rem*}
  Take indices $i_0,i_1,i_2,j_1,j_2,k_1,k_2,l_1,l_2,m_1,m_2 \in \{1
  \upto n\}$ as in the hypothesis of \tref{test}. Explicit computation
  shows that
  \[
  g\left(d_{\{i_0,i_1\}},d_{\{j_1,j_2\} },d_{\{ k_1,k_2\}},
    d_{\{l_1,l_2\}},d_{\{m_1,m_2\}},d_{\{i_0,i_2\}}\right),
  \]
  viewed as a polynomial in variables $d_{\{i,j\}}$, contains the term
  $2 d_{\{i_0,i_1\}}^2 d_{\{i_0,i_2\}}$. Notice that the index $i_0$
  occurs three times in this term (when writing it out as a product
  rather than squaring the first variable). It follows from
  \mycite[Proposition~2.2(b) and Lemma~2.3]{FirstPaper} that this term
  does not occur in any relationship of degree~3 between the
  $d_{\{i,j\}}$. In particular, $g\left(d_{\{i_0,i_1\}},
    d_{\{j_1,j_2\}},d_{\{ k_1,k_2\}},d_{\{l_1,l_2\}},d_{\{m_1,m_2\}},
    d_{\{i_0,i_2\}}\right)$ is not a relationship between the
  $d_{\{i,j\}}$. It follows that there exists a dense, open subset
  $\Omega \subseteq \left(\RR^2\right)^n$ such that for all point
  configurations $(p_1 \upto p_n) \in \Omega$ the hypotheses of
  \tref{test} are met. This provides a new proof for the fact that
  ``most'' point configurations are reconstructible from distances,
  which appeared in greater generality in~[\citenumber{FirstPaper},
  Theorem~2.6].
\end{rem*}

How many tests do we have to conduct for checking that the conditions
in~\eqref{eqTest} are satisfied? There are~$n$ choices for $i_0$, the
index that is repeated. For each choice of $i_0$, there are
$(n-1)(n-2)$ choices for $i_1$ and $i_2$ (since these three indices
must be distinct). Having chosen $i_0$, $i_1$, and $i_3$, there are
$\binom{n}{2} - 2$ choices for the set $\{j_1,j_2\}$, $\binom{n}{2} -
4$ choices for the set $\{k_1,k_2\}$ and so on. Altogether, we obtain
\begin{multline*}
  n (n-1) (n-2) \left(\binom{n}{2} - 2\right) \left(\binom{n}{2} -
    3\right) \left(\binom{n}{2} - 4\right) \left(\binom{n}{2} -
    5\right) = \\
  \frac{1}{16} \left(n^{11}-7 n^{10}-8 n^9+138 n^8-83 n^7-983 n^6+1074
    n^5+2996 n^4-3672 n^3-3296 n^2+ 3840 n\right)
\end{multline*}
choices.

\begin{cor}
  There exists an open and dense set $\Omega \subset (\RR^2)^n$ of
  reconstructible $n$-point configurations and an algorithm in
  $O(n^{11})$ steps to determine whether any $( p_1,\ldots,p_n )\in
  (\RR^2)^n $ lies in $\Omega$.
\end{cor}

\begin{remark}
\label{algo_dim_m}
The algorithm given by \tref{test} can be
generalized to $\RR^m$ if $n \geq m+2$. For each choice of $m+2$
indices $i_0 \upto i_{m+1}$ we have the relationship
\[
\det\left(d_{\{i_\nu,i_\mu\}}-d_{\{i_\nu,i_0\}} -
  d_{\{i_\mu,i_0\}}\right)_{\nu,\mu = 1 \upto m+1} = 0,
\]
which can be expressed as $g_m\left(d_{\{i_0,i_1\}} \upto
  d_{\{i_m,i_{m+1}\}}\right) = 0$ with $g_m$ an appropriate polynomial
in $k := \binom{m+2}{2}$ variables. Now we obtain a generalization of
\tref{test} which says that if for all pairwise distinct choices $S_1
\upto S_k \in {\mathcal P}$ with $S_1 \cap S_k \ne \emptyset$ we have
\begin{equation} \label{eqGen}
  g_m\left(d_{S_1} \upto d_{S_k}\right) \ne 0,
\end{equation}
then the configuration $p_1 \upto p_n$ is reconstructible from
distances. We see that there are
\[
n (n-1) (n-2) \prod_{j=2}^{k-1} \left(\binom{n}{2} - j\right) =
O\left(n^{m^2+3m+1}\right)
\]
steps for checking the reconstructibility of $p_1 \upto p_n$. It also
follows from \mycite[Proposition~2.2(b) and Lemma~2.3]{FirstPaper}
that there exists a dense open subset $\Omega \subseteq
\left(\RR^m\right)^n$ where the inequalities~\eqref{eqGen} are all
satisfied.
\end{remark}

\section{Numerical Experiments}
A simple Matlab code (see the appendix)
was used to check for the reconstructibility 
of some $n$-point configurations.
In the code, we traded simplicity for speed
in an attempt to make the algorithm more easily understandable.
Even so, 
we were able to show that some $n$-point configurations were reconstructible,
with $n=5,6,7$ and even $8$
in a reasonable time.
Corresponding CPU times
and number of combinations to be checked 
are given in Table \ref{TableTime}. 
The computations were done using Matlab version 6.1 on a 
Sun (4$\times$ultraSPARC-II, 480 MHZ).

\begin{table}
\centering          
\begin{tabular}{|crr|}      
\hline 
$n$ & \#  combinations & CPU time in seconds   \\    
\hline  
5             &        100,800   &     72   \\ 
6             &      2,059,200   &  1,170   \\
7             &     19,535,040   &  9,920 \\
8             &    120,556,800   & 58,375 \\
 \hline
\end{tabular}
\label{TableTime}          
\caption{Time required to check for the reconstructibility of an $n$-point configuration.} 
\end{table}

An important point to observe is that 
if a point configuration fails to satisfy
one of the conditions in (\ref{eqTest}),
it does not mean that it is not reconstructible.
For example, it is not hard to show that every square is
reconstructible (see \mycite[Example~2.12]{FirstPaper}). 
But, as one can check, 
squares satisfy neither (\ref{eqTest}) nor (\ref{eqTest2}).
This is due to the fact that squares have repeated distances.
Indeed, any planar $n$-point configuration with repeated distances
will fail the reconstructibility test.
(See \mycite{FirstPaper} for a proof of this fact 
and ideas on how to modify the algorithm 
to take care of point configurations with repeated distances.)
Also, the point configuration given by
\[
p_1=(0,0), p_2=(7,0), p_3=(5,-1), p_4=(3,-3), p_5=(11,2)
\]
does not satisfy  (\ref{eqTest}), 
even though its pairwise distances are all distinct.  
However, one can show that it is actually reconstructible. 
(It suffices to show that the permutations of the distances 
which make $g$ equal to zero all violate one of the relationships 
that exist between the pairwise distances of five points. 
We checked this numerically.)
Our test is thus not perfect.

Observe that,
when using points with small integer coordinates, 
the polynomial $g$ can be evaluated exactly on a computer.
We can thus determine precisely whether such a 
point configuration satisfies the conditions of (\ref{eqTest}).
An interesting question is:
given a planar $n$-point configuration with integer coordinates
and lying inside the box $[0,N]\times [0,N]$,
what are the chances that it will fail the reconstructibility test?
Numerical experiments showed that it is quite likely, 
even when configurations with repeated distances are excluded.
For $N=3$, 
we found that about 61\% of configurations
of four points whose distances are not repeated fail the test.
(More precisely, we generated all possible 
$p_1=(x_1,y_1)$, $p_2=(x_2,y_2)$, $p_3=(x_3,y_3)$, $p_4=(x_4,y_4)$ 
with coordinates in $\{ 0,1,2,3\}$ and such that
either $x_i< x_{i+1}$ or $x_i=x_{i+1}$ and $y_i<y_{i+1}$, for all $i=1,2,3,4$.
Of those 1820 four-point configurations,
we found that 1636 had repeated distances 
while a total of 1748 failed the test.)
For $N=4$, this percentage went down to about 30\%, which is still quite high.

It would be interesting 
to determine whether such high rates of failure 
are also 
observed when the coordinates of the points are not necessarily integers. 
But, in general, floating-point arithmetic prevents us for 
determining whether a polynomial function is exactly zero.
We must thus replace the $g=0$ in conditions 1 and 2
by $|g |\leq \epsilon$, 
for some $\epsilon$ determined by the machine precision and possible
noise in the measurements.
However, numerical tests have shown that 
if the coordinates of four points are chosen randomly in $(0,1)$ 
(using the Matlab {\em rand} function), 
then the polynomial $g$ in (\ref{eqTest}) rarely takes very small values.
For example, after generating 5000 different random four-point configurations, 
we found that only 22 of those generated a $g$ 
with a value less than $10^{-7}$. 
In another set of 5000 four-point configurations, 
we found only 6 which generated a $g$ with a value less than $10^{-8}$.
In a final set of 10,000 four-point configurations,
we found none which generated a $g$ with a value less than $10^{-9}$.
As these values are well above the maximal error expected 
with such data when evaluating $g$ using Matlab,
this implies that none of the 20,000 random four-point configurations
we generated could possibly fail the test.


\section{The Case of Orientation Preserving Rigid Motions in the Plane}
In the previous two sections, we considered the case
where the shape of an n-point configurations
is defined by $p_1,\ldots,p_n\in\RR^m$ up to rigid motions.
Recall that the group of rigid motions in $\RR^m$, 
sometimes called the Euclidean group and denoted by $E(m)$, 
is generated by rotations, translations and reflections in $\RR^m$.
However, in certain circumstances, 
it may be desirable to be able to determine 
whether two point configurations 
are equivalent up to strictly orientation preserving rigid motions.
The group of orientation preserving rigid motions,
sometimes called the {\em special Euclidean group} and denoted by $SE(m)$,
is the one that is generated by rotations and translations in $\RR^m$.

For simplicity, we again restrict ourselves to the planar case $m=2$.
Given a planar point configuration $p_1,\ldots,p_n\in \RR^2$,
we would like to be able to determine whether any other 
planar $n$-point configuration $q_1,\ldots,q_n$
is the same as $p_1,\ldots,p_n$ up to a rotation and a translation? 
Given any $q_i,q_j,q_k$ in the plane,
denote by $a_{q_i,q_j,q_k}$
the signed area of the  parallelogram spanned by $q_i-q_j$ and
$q_k-q_j$, so
\[
a_{q_i,q_j,q_k} = \det(q_i-q_k,q_j-q_k).
\]
Since signed areas are unchanged under rotations and translations,
the function $I:\RR^2\times\RR^2\times\RR^2\times\RR^2\rightarrow \RR$
defined by
\begin{multline}
I(q_1,q_2,q_3,q_4)=(a_{q_1,q_2,q_4}^2-a_{q_1,q_3,q_4}^2)
(a_{q_1,q_2,q_3}^2-a_{q_1,q_3,q_4}^2)\\
(a_{q_1,q_2,q_3}^2-a_{q_1,q_2,q_4}^2)
(a_{q_1,q_2,q_3}-a_{q_1,q_2,q_4}+ 2a_{q_1,q_3,q_4})\\
(a_{q_1,q_2,q_3}-2 a_{q_1,q_2,q_4}+ a_{q_1,q_3,q_4})
(2 a_{q_1,q_2,q_3}- a_{q_1,q_2,q_4}+ a_{q_1,q_3,q_4})
\end{multline}
is invariant under the action of $SE(2)$.
Moreover, one can check that
it is also invariant under a relabeling of the four points
$q_1,q_2,q_3,q_4$.
However, it is {\em not} invariant 
under rigid motions in general.
Indeed, any transformation which is a rigid motion 
but does {\em not} preserve the orientation 
will transform $I$ into $-I$.

Given an $n$-point configuration $q_1,\ldots,q_n$ with $n\geq 4$,
we can evaluate $I$ on all possible subsets of four points of 
$\{q_1,\ldots,q_n \}$. 
We consider the distribution of the value of these $I$'s, 
i.e.~the distribution of the 
\[
I_{i_1,i_2,i_3,i_4}=I(q_{i_1},q_{i_2},q_{i_3},q_{i_4}),\text{ for all }
i_1<i_2<i_3<i_4\in\{1,\ldots,n \}.
\]

\begin{prop}
Let $n\geq 4$
and let $p_1,\ldots,p_n \in \RR^2$ be an $n$-point configuration 
which is reconstructible from distances.
Assume that the distribution of the $I$'s of this point configuration
is not a symmetric function
(i.e.~that the distribution of the $I$'s is not the same 
as the distribution of the $-I$'s.)
Let $q_1,\ldots,q_n\in \RR^2$ be another $n$-point configuration.
Then both the distribution of the distances and
the distribution of the $I$'s
of the two point configurations are the same
if and only if
there exists a rotation and a translation which maps
one point configuration onto the other.
\end{prop}
\begin{proof}
Observe that, 
in addition to being invariant under rotations and translations of the points, 
the distribution of the value of the $I$'s
is also independent of the labeling of the points.  
The same holds for the distribution of pairwise distances.
So if two $n$-point configurations are the same up to a rotation, 
a translation and a relabeling, 
then the distribution of the $I$'s and the distribution of the distances 
are the same for both.
Thus the {\em if} is clear.

Now assume that the distribution of the distances 
and the distribution of the $I$'s 
are the same for both point configuration.
Since $p_1,\ldots,p_n$ is, by hypothesis, reconstructible, 
this implies that there exists a rigid motion $(M,T)$ 
and a relabeling $\pi:\{1,\ldots,n \} \hookrightarrow \{1,\ldots,n \}$
such that $M p_i+T=q_{\pi(i)}$, for all $i=1,\ldots,n$.
If $(M,T)$ is not in $SE(2)$, then it maps each 
$I(p_{i_1},p_{i_2},p_{i_3},p_{i_4})$ to $-I(p_{i_1},p_{i_2},p_{i_3},p_{i_4}))$.
But this is a contradiction, 
since the distribution of the $I's$ is not symmetric.
Thus $g$ is in $SE(2)$.
This shows the {\em only if}. 
\end{proof}

\begin{remark}
One can actually show that if  
$p_1,\ldots,p_4$ is equivalent to $q_1,\ldots,q_4$ up to a rigid motion,
then $p_1,\ldots,p_4$ is equivalent to $q_1,\ldots,q_4$ 
up to a rotation and a translation if and only if 
$I(p_1,\ldots,p_4)=I(q_1,\ldots,q_4)$.
(Indeed, $I$ is one of the two fundamental invariants
of the action of $SE(2)\times S_4$ 
on $\RR^2\times\RR^2\times\RR^2\times\RR^2$
which we obtained 
using the invariant theory package in Magma~[\citenumber{Magma}].
By construction, 
these two invariants thus distinguish the orbits of $SE(2)\times S_3$.
The other invariant is actually 
unchanged under the action of the full Euclidean group 
$E(2)$ and so $I$ alone distinguishes the orbits of $SE(2)$
within the orbits of $E(2)$. 
\end{remark}

\section{The Case of Rotations, Translations and Scalings}

In certain circumstances, 
it may also be desirable 
to be able to determine whether two point configurations
are the same up to a rigid motion and a scaling.
This can be done using a simple variation of the previous approach.
Given a distribution of distances $\{ d_{\{ i,j \}}=\|p_i-p_j \|^2\}$,
let $d_{max}$ be the largest distance
\[
d_{max}=\max \{ d_{\{ i,j \}} | \{i,j \}\in {\mathcal P}\},
\]
which can be assumed to be non-zero since otherwise all points
coincide. We can consider the distribution of the rescaled distances 
$\{ \frac{d_{\{ i,j\}}}{d_{max}}\}_{\{i,j \}\in {\mathcal P}}$.
In addition to being invariant under rigid motions and relabeling,
the distribution of the rescaled distances
is also invariant under a scaling of the points 
\[
p_i\mapsto \lambda p_i,\text{ for every }i=1,\ldots,n.
\]
for any $\lambda \in \RR_{\neq 0}$.

\begin{prop}
Let $n\geq m+2$.
There exists an open, dense subset $\Omega$ of $(\RR^m)^n$ such that
if an $n$-point configuration $p_1,\ldots,p_n$
is such that $(p_1,\ldots,p_n)\in\Omega$,
then $p_1,\ldots,p_n$ is uniquely determined, 
up to rotations, translations, reflections, scalings and relabeling of the points, 
by the distribution of its rescaled pairwise distances 
$\{ \frac{ d_{\{ i,j \}}}{d_{max}} \}_{\{i,j\}\in {\mathcal P}}$.
Moreover, there is an algorithm in $O(n^{\frac{m^2+3m+12}{2}})$
steps to determine whether $(p_1,\ldots,p_n)\in\Omega$.
\end{prop}

\begin{proof}
Let $p_1,\ldots,p_n\in\RR^m$ 
be an $n$-point configuration which is 
reconstructible from distances
and whose pairwise distances are not all zero.
Observe that if  $q_1,\ldots,q_n \in \RR^2$ is another $n$-point configuration,
then the distributions of the rescaled distances
of both point configurations are the same
if and only if 
there exists a rigid motion followed by a scaling which maps
one point configuration onto the other.
The claim is thus a direct corollary of Theorem~2.6 from \mycite{FirstPaper}
and of Remark \ref{algo_dim_m}.
\end{proof}

\section*{Acknowledgments}

Mireille Boutin wishes to thank David Cooper and Senem
Velipasalar~[\citenumber{TasdizenVelipasalarCooper}] who gave her the
initial motivation for this research, Jean-Philippe Tarel for his
constructive remarks on this approach and her co-author for being such
a joy to work with. Gregor Kemper expresses his gratitude to his hosts
at Purdue University, where part of this research took place.  Another
part of the research was carried out while both authors were visiting
the Mathematical Sciences Research Institute in Berkeley.  We wish to
thank Michael Singer and Bernd Sturmfels for the invitation.


\section*{Appendix}

\% This is a simple Matlab function that determines whether\\
\% the pairwise distances between the points of the planar $n$-point configuration \\
\% defined by the $n$ columns of a $2$-by-$n$ matrix $p$ satisfy the conditions of (\ref{eqTest})\\
\% with $0$ replaced by some (small) number e. \\

\noindent function  g=evaluate\_g(p,e) \\
\% p is a $2$-by-$n$ matrix. \\
\% e should be chosen depending on the machine precision.\\
\% This function returns $0$ if $g\leq e$ for all sets of pairs \\
\% of the conditions in \ref{eqTest} and 1 otherwise.

\noindent n=length(p);\\
\noindent\% Compute the squares of the distances d.\\
\noindent for i=1:n-1, for j=i+1:n \\
\phantom{d} d(i,j)=sum((p(i,:)-p(j,:)).\^{}2);\\
end, end \\
d(n,n)=0; \\
d=d+d';\\
for i0=1:n, for i1=1:n, if i1\~{}=i0 \\
\phantom{d} x12=d(i0,i1);\\
\phantom{d} for i2=1:n, if i2\~{}=i1 \& i2\~{}=i0\\
\phantom{dd}  x34=d(i0,i2); \\   
\phantom{dd}  for j1=1:n-1, for j2=j1+1:n \\
 \phantom{dd} if (j1\~{}=i0 $|$ j2\~{}=i1) \& (j1\~{}=i1 $|$ j2\~{}=i0) \& (j1\~{}=i0 $|$ j2\~{}=i2) \& (j1\~{}=i2 $|$ j2\~{}=i0) \\
\phantom{ddd}  x13=d(j1,j2);\\
\phantom{ddd}   for k1=1:n-1, for k2=k1+1:n;\\
\phantom{ddd}   if (k1\~{}=i0 $|$ k2\~{}=i1) \& (k1\~{}=i1 $|$ k2\~{}=i0) \& 
(k1\~{}=i0 $|$ k2\~{}=i2) \& (k1\~{}=i2 $|$ k2\~{}=i0) \& \ldots \\
\phantom{ddd} (k1\~{}=j1 $|$ k2\~{}=j2) \& (k1\~{}=j2 $|$ k2\~{}=j1) \\
\phantom{dddd}    x14=d(k1,k2);\\
\phantom{dddd}    for l1=1:n-1, for l2=l1+1:n \\ 
\phantom{dddd}    if (l1\~{}=i0 $|$ l2\~{}=i1) \& (l1\~{}=i1 $|$ l2\~{}=i0) \& (l1\~{}=i0 $|$ l2\~{}=i2) \& (l1\~{}=i2 $|$ l2\~{}=i0) \& \ldots \\
\phantom{dddd}    
 (l1\~{}=j1 $|$ l2\~{}=j2) \& (l1\~{}=j2 $|$ l2\~{}=j1) \&
 (l1\~{}=k1 $|$ l2\~{}=k2) \& (l1\~{}=k2 $|$ l2\~{}=k1) \\
\phantom{ddddd}     x23=d(l1,l2);\\
\phantom{ddddd}       for m1=1:n-1, for m2=m1+1:n \\
\phantom{ddddd}     if (m1\~{}=i0 $|$ m2\~{}=i1) \& (m1\~{}=i1 $|$ m2\~{}=i0) \& (m1\~{}=i0 $|$ m2\~{}=i2) \& \ldots \\
\phantom{ddddd}    (m1\~{}=i2 $|$ m2\~{}=i0) \&  
 (m1\~{}=j1 $|$ m2\~{}=j2) \& (m1\~{}=j2 $|$ m2\~{}=j1) \&\ldots \\
\phantom{ddddd} (m1\~{}=k1 $|$ m2\~{}=k2) \& (m1\~{}=k2 $|$ m2\~{}=k1) \&
(m1\~{}=l1 $|$ m2\~{}=l2) \& (m1\~{}=l2 $|$ m2\~{}=l1)\\
\phantom{dddddd} x24=d(m1,m2); \\
\phantom{dddddd} m11=-2*x14; \\
\phantom{dddddd} m12=x12-x14-x24; \\
\phantom{dddddd} m13=x13-x14-x34; \\
\phantom{dddddd} m22=-2*x24; \\
\phantom{dddddd} m23=x23-x24-x34; \\
\phantom{dddddd} m33=-2*x34; \\
\phantom{dddddd} Mu=
m11*m22*m33-m11*m23\^{}2-m12\^{}2*m33-m12*m23*m13+m13*m12*m23-m22*m13\^{}2;\\
\phantom{dddddd}      if Mu== 0 \\ 
 \phantom{ddddddd}      g=0; \\
\phantom{ddddddd}       return \\
\phantom{dddddd}      end, end, end, end, end, end, end, end, end, end, end, end, \\
end, end, end, end, end, end \\
g=1;

\addcontentsline{toc}{section}{References}
\bibliographystyle{mybibstyle} 
\bibliography{ourbib}

\bigskip

\begin{center}
\begin{tabular}{lll}
  Mireille Boutin & &Gregor Kemper \\
  Department of Mathematics  & & Technische Universit\"at M\"unchen \\
  Purdue University & & Zentrum Mathematik - M11 \\
  150 N. University St.  & & Boltzmannstr. 3 \\
  West Lafayette, IN & & 85\,748 Garching \\
  47907, USA & & Germany \\
  {\tt boutin$@$math.purdue.edu} & & {\tt kemper$@$ma.tum.de}
\end{tabular}
\end{center}

\end{document}